# Analytic properties of hidden variable recurrent fractal interpolation function with function contractivity factors

Mi-Kyong Ri, Chol-Hui Yun

Faculty of Mathematics, **Kim Il Sung** University,

Pyongyang, Democratic People's Republic of Korea

Abstract: In this paper, we analyze the smoothness and stability of hidden variable recurrent fractal interpolation functions (HVRFIF) with function contractivity factors introduced in Ref. 1. The HVRFIF is a hidden variable fractal interpolation function (HVFIF) constructed by recurrent iterated function system (RIFS) with function contractivity factors. An attractor of RIFS has a local self-similar or self-affine structure and looks more naturally than one of IFS. The contractivity factors of IFS(RIFS) determine fractal characteristic and shape of its attractor. Therefore, the HVRFIF with function contractivity factors has more flexibility and diversity than the HVFIF constructed by iterated function system (IFS) with constant contractivity factors. The analytic properties of the interpolation functions play very important roles in determining whether these functons can be applied to the practical problems or not. We analyze the smoothness of the one variable HVRFIFs in Ref. 1 and prove their stability according to perturbation of the interpolation dataset.

## 1. Introduction

A fractal interpolation function (FIF) is an interpolation function whose graph is a fractal set. In 1986, M.F.Barnsley defined a fractal interpolation function as an interpolation function whose graph is an attractor of an iterated function system (IFS) in Ref. 2. The FIFs can model better the irregular and complicated phenomena in nature than the classical interpolation functions such as polynomials and spline. Therefore, construction[3, 4, 5], fractal dimension[6, 3, 4] and analytic properties[6, 7, 8, 9] of FIFs have been widely studied and applied to a lot of fields such as approximation theory, computer graphics, image compression, metallurgy, physics, geography, geology, data fitting, signal process and so on. The contractivity factors of IFS determine fractal characteristic and shape of its attractor. Therefore, the constructions of FIFs with constant contractivity factors lack the flexibility and have been developed into one with function contractivity factors. Many researchers have studied construction,[10, 11, 1] dimension[3, 12, 13] and analytic properties[14, 15, 16, 17, 5] of the FIFs with function contractivity factors.

A recurrent fractal interpolation function (RFIF) is an interpolation function whose graph is an attractor of a recurrent iterated functions system (RIFS) defined by M.F.Barnsley *et al.* in Ref. 18. The RIFS is a generalization of IFS and constructs local self-similar sets which have more complicated structures than self-similar sets constructed by IFSs. A.E.Jacquin[19] presented a new method of fractal image compression (FIC) based on RIFSs, which is a more flexible than one using IFSs. P. Bouboulis *et al.*[20] introduced a method of using RFIFs to improve the FIC. Metzler and Yun[21] presented a new construction of RIFSs with function contractivity factors instead of constant contractivity factors in earlier method[20, 19] and applied the RIFSs to image compression to improve the quality of decoded image. P. Bouboulis *et al.*[22] introduced a construction of recurrent bivariate fractal interpolation functions and computed their box-counting dimension. Yun *et al.*[12] constructed the RFIFs using the RIFS with function contractivity factors and estimated the box-counting dimension of the graph of the constructed RFIF. They also constructed bivariate RFIFs with function contractivity factors and estimated the box-counting dimension of the graph of the RFIFs[16]. D.C.Luor[23] introduced a construction of FIFs with locally self-similar graphs in $\mathrm{R}^2$, where contractivity factors are functions and homeomorphisms, that map domains to regions, are not contractions.

Barnsley *et al.*[24] introduced a concept of hidden variable fractal interpolation function (HVFIF) whose graph is not self-similar and self-affine set. We can use the hidden variables to control fractal characteristics and shapes of the graphs of HVFIFs more flexibly. Study on the analytic properties of HVFIFs guarantees possibility of applying to modeling the practical problems. The construction and analytic properties of HVFIFs with constant contractivity factors were studied in Refs. 25, 26, 27, 28 and 17. To generalize the method, researchers have studied construction,[29, 11] smoothness and stability[17] of HVFIFs with function contractivity factors. In order to ensure the flexibility and diversity of construction of FIFs, construction of one variable and bivariable HVRFIFs using RIFS with four function contractivity factors was presented in Ref. 1.

We analyze the smoothness and stability of one variable HVRFIFs constructed in Ref. 1. This paper is organized as follows:

In section 2, we introduce construction of HVRFIFs in Ref. 1 (Theorem 1, 2). In section 3, we analyze the smoothness of the constructed HVRFIFs and in section 4, prove their stability (Theorem 4-7). These results become theoretical basis of practical application such as approximation theory, computer graphics and image process.

## 2. Hidden Variable recurrent fractal interpolation function(HVRFIF)

In this section, we introduce the construction of one variable HVRFIFs presented in Ref. 1 which is used to study their analytic properties.

Let a data set $P_0$ in $\mathrm{R}^2$ be given by

$$P_0 = \{(x_i, y_j) \in \mathrm{R}^2; i = 0, 1, \cdots, n\}, \quad (-\infty < x_0 < x_1 < \cdots < x_n < +\infty)\,. \tag{1}$$

We extend the dataset as follows:

$$P = \{(x_i, y_i, z_i,) = (x_i, \vec{y}_i) \in \mathrm{R}^3; i = 0, 1, \cdots, n\},\ (-\infty < x_0 < x_1 < \cdots < x_n < \infty), \tag{2}$$

where $\vec{y}_i = (y_i, z_i)$ and $g_i'$, $i = 0, 1, \cdots, n$ are parameters. Moreover, we denote $I_i = [x_{i-1}, x_i]$ and $I = [x_0, x_n]$ where $I_i$, $i \in \{1, \cdots, n\}$ is called a region. Let $l$ be an integer with $2 \le l \le n$. We make subintervals $\tilde{I}_k$, $k = 1, \cdots, l$ of $I$ consisting of some regions and call $\tilde{I}_k$ domain. Let we denote start point and end point of $\tilde{I}_k$ by $s(k)$, $e(k)$, respectively. Then we get the following mappings:

$$s: \{1, \cdots, l\} \to \{1, \cdots, n\},\ e: \{1, \cdots, l\} \to \{1, \cdots, n\},$$

and $\tilde{I}_k$ is denoted by $\tilde{I}_k = [x_{s(k)}, x_{e(k)}]$. We suppose that $e(k) - s(k) \ge 2$, $k = 1, \cdots, l$ which means that the interval $\tilde{I}_k$ contains at least 2 of $I_i$ s.

For each $i\,(1 \le i \le n)$, we take a $k (\in \{1, \cdots, l\})$ and denote it by $\gamma(i)$. Let mappings $L_{i,k}: [x_{s(k)}, x_{e(k)}] \to [x_{i-1}, x_i]$ be contraction homeomorphisms that map end points of $\tilde{I}_k$ to end points of $I_i$, i.e. $L_{i,k}(\{x_{s(k)}, x_{e(k)}\}) = \{x_{i-1}, x_i\}$.

We define mappings $F_{i,k}:\tilde{I}_k\times R^2\to R,\ i=1,\cdots,n$ as follows:

$$\vec{F}_{i,k}(x,\ \vec{y})=\begin{pmatrix} s_i(L_{i,k}(x))y+s_i'(L_{i,k}(x))z+q_{i,k}(x) \\ \tilde{s}_i(L_{i,k}(x))y+\tilde{s}_i'(L_{i,k}(x))z+\tilde{q}_{i,k}(x)\end{pmatrix}=\begin{pmatrix} s_i(L_{i,k}(x)) & s_i'(L_{i,k}(x)) \\ \tilde{s}_i(L_{i,k}(x)) & \tilde{s}_i'(L_{i,k}(x))\end{pmatrix}\begin{pmatrix} y \\ z\end{pmatrix}+\begin{pmatrix} q_{i,k}(x) \\ \tilde{q}_{i,k}(x)\end{pmatrix},$$

where $s_i,\ s_i',\ \tilde{s}_i,\ \tilde{s}_i':I_i\to R$ are Lipshitz functions on $I_i$ whose absolute value is less than 1 (which are called contractivity factors) and $q_{i,k}$ , $\tilde{q}_{i,k}:\tilde{I}_k\to R$ are Lipshitz functions such that if $\alpha\in\{s(k),e(k)\}$ , $L_{i,k}(x_\alpha)=x_a$ $a\in\{i-1,i\}$, then $\vec{F}_{i,k}(x_\alpha,\ \vec{y}_\alpha)=\vec{y}_a$. Then, $\vec{F}_{i,k}(x,\ \vec{z})$ is obviously Lipshitz function.

Example 1. An example of $q_{i,k}$ , $\tilde{q}_{i,k}$ satisfied the above conditions is as follows:

$$q_{i,k}(x)=-s_i(L_{i,k}(x))g_{i,k}(x)-s_i'(L_{i,k}(x))g_{i,k}'(x)+h_i(L_{i,k}(x)),$$

$$\tilde{q}_{i,k}(x)=-\tilde{s}_i(L_{i,k}(x))g_{i,k}(x)-\tilde{s}_i'(L_{i,k}(x))g_{i,k}'(x)+h_i(L_{i,k}(x)),$$

$$g_{i,k}(x)=\frac{x-x_{s(k)}}{x_{e(k)}-x_{s(k)}}y_{e(k)}+\frac{x-x_{e(k)}}{x_{s(k)}-x_{e(k)}}y_{s(k)},\quad g_{i,k}'(x)=\frac{x-x_{s(k)}}{x_{e(k)}-x_{s(k)}}z_{e(k)}+\frac{x-x_{e(k)}}{x_{s(k)}-x_{e(k)}}z_{s(k)},$$

$$h_i(x)=\frac{x-x_{i-1}}{x_i-x_{i-1}}y_i+\frac{x-x_i}{x_{i-1}-x_i}y_{i-1}.$$

Let $D\subset R^2$ be a sufficiently large bounded set containing $\vec{y}_i,\ i=1,\cdots,n$.

We define transformations $\vec{W}_i:\tilde{I}_{\gamma(i)}\times D\to I_i\times R^2,\ i=1,\cdots,n$ by

$$\vec{W}_i(x,\vec{y})=(L_{i,\gamma(i)}(x),\ \vec{F}_{i,\gamma(i)}(x,\vec{y})),\ i=1,\cdots,n.$$

Then, as we know from the definitions of $L_{i,\gamma(i)}$ and $\vec{F}_{i,\gamma(i)}$, $\vec{W}_i$ maps the data points on end of the domain $\tilde{I}_{\gamma(i)}$ to the data points of $I_i$. For a function $f$, let us denote $\bar{f}=\max_x|f(x)|$. We denote $\bar{S}=\max\{\bar{s}_i+\bar{\tilde{s}}_i,\ \bar{s}_i'+\bar{\tilde{s}}_i';\ i=1,\cdots,n\}$.

**Theorem 1.**[1] **If $\bar{S}<1$, then there exists some distance $\rho_\theta$ equivalent to the Euclidean metric such that $\vec{W}_i(i=1,\cdots,n)$ are contraction transformations with respect to the distance $\rho_\theta$.**

We define a row-stochastic matrix $M=(p_{st})_{n\times n}$ by

$$p_{st}=\begin{cases}1/a_s, & I_s\subseteq\tilde{I}_{\gamma(t)} \\ 0, & I_s\not\subset\tilde{I}_{\gamma(t)}\end{cases},$$

where for every $s(1\le s\le n)$, the number $a_s$ indicates the number of the domains $\tilde{I}_k$ containing the region $I_s$, which means that $p_{st}$ is positive if there is a transformation $W_i$ mapping $I_s$ to $I_t$.

Then, we have RIFS$\{R^3;\ M;\ W_i,i=1,\cdots,n\}$ corresponding to the extended dataset *P*.

### 2.2 Construction of HVFIF

For the RIFS constructed above, we have the following theorem.

**Theorem 2.**[1] **There is a continuous function $\vec{f}$ interpolating the extended data set $P$ such that the graph of $\vec{f}$ is the attractor A of RIFS constructed above.**

In the theorem 2, the following function space is considered and the Read-Bajraktarevic operator is as follows:

$$\bar{C}(I)=\{\vec{h}:I\to R^2;\vec{h}\ \text{interpolates the extended data set } P \text{ and is continuous}\}$$

$$(T\vec{h})(x)=\vec{F}_{i,\gamma(i)}(L_{i,\gamma(i)}^{-1}(x),\ \vec{h}(L_{i,\gamma(i)}^{-1}(x))),\quad x\in I_i.$$

We denote a unique fixed point of the Read-Bajraktarevic operator by

$$\vec{f}(x)=\vec{F}_{i,\gamma(i)}(L_{i,\gamma(i)}^{-1}(x),\ \vec{f}(L_{i,\gamma(i)}^{-1}(x))).$$

Let us denote the vector valued function $\vec{f}:I\to R^2$ by $\vec{f}=(f_1,\ f_2)$, where $f_1:I\to R$ interpolates the given dataset $P_0$, which is called a hidden variable recurrent fractal interpolation function (HVRFIF). Furthermore, a set $\{(x,\ f_1(x)):x\in I\}$ is a projection of $A$ on $R^2$. Then we have $\vec{f}(x,\ y)=\vec{F}_{i,\gamma(i)}(L_{i,\gamma(i)}^{-1}(x),\ \vec{f}(L_{i,\gamma(i)}^{-1}(x))),\quad x\in I_i$, i.e.

$$\vec{f}(x)=\vec{F}_{i,\gamma(i)}(L_{i,\gamma(i)}^{-1}(x),\ f_1(L_{i,\gamma(i)}^{-1}(x)),\ f_2(L_{i,\gamma(i)}^{-1}(x))),\quad x\in I_i.$$

Therefore, for all $x\in I$, the HVRFIF $f_1$ satisfies (3) and $f_2$ satisfies (4)

$$f_1(x)=s_i(x)f_1(L_{i,\gamma(i)}^{-1}(x))+s_i'(x)f_2(L_{i,\gamma(i)}^{-1}(x))+q_{i,\gamma(i)}(L_{i,\gamma(i)}^{-1}(x)),\ x\in I_i \tag{3}$$

$$f_2(x)=\tilde{s}_i(x)f_1(L_{i,\gamma(i)}^{-1}(x))+\tilde{s}_i'(x)f_2(L_{i,\gamma(i)}^{-1}(x))+\tilde{q}_{i,\gamma(i)}(L_{i,\gamma(i)}^{-1}(x)),\ x\in I_i. \tag{4}$$

### 3. Smoothness of the HVRFIF.

In this section, we analyze a smoothness of the constructed HVRFIF, i.e. Hoelder index of the HVRFIF is calculated in theorem 3.

Let $I=[0,1]$ and $L_{i,k}:\tilde{I}_k\to I_i$ be similitudes. Let us denote as follows: $\bar{g}=\max_x|g(x)|$, $\bar{s}=\max_i\{\bar{s}_i\}$, $\bar{\tilde{s}}=\max_i\{\bar{\tilde{s}}_i\}$, $\bar{s}'=\max_i\{\bar{s}_i'\}$, $\bar{\tilde{s}}'=\max_i\{\bar{\tilde{s}}_i'\}$, $\omega_k=\max\{\bar{s}_k,\ \bar{s}_k'\}$, $\tilde{\omega}_k=\max\{\bar{\tilde{s}}_k,\ \bar{\tilde{s}}_k'\}$,

$$\mathrm{L}_L=\max_{k=1,n}\left\{\frac{|I_k|}{|\tilde{I}_{\gamma(k)}|}\right\},\quad l_L=\min_{k=1,n}\left\{\frac{|I_k|}{|\tilde{I}_{\gamma(k)}|}\right\},\quad l_{\tilde{I}}=\min_k\{|\tilde{I}_k|\}.$$

**Lemma 1. If $0<\alpha$, then for $x$ obeying $0<x<1$,**

$$0<-x^\alpha\ln x\le\frac{1}{\alpha e}.$$

**Proof.** Since $0<x<1$, we have $0<-x^\alpha\ln x$. Let us define a function $f:\ [0,1]\to R$ by $f(x)=x^\alpha\ln x+\frac{1}{\alpha e}$. Then $f'(x)=\alpha x^{\alpha-1}\ln x+x^{\alpha-1}=x^{\alpha-1}(\alpha\ln x+1)$. Therefore, if $x>e^{-\frac{1}{\alpha}}$, then $f'(x)>0$ and if $0<x<e^{-\frac{1}{\alpha}}$, then $f'(x)<0$. Hence, for all $x$ obeying $0<x<1$, we have $f(x)\ge f(e^{-\frac{1}{\alpha}})=0$, i.e. $x^\alpha\ln x+\frac{1}{\alpha e}\ge 0$.□

**Theorem 3. Let $\vec{f}=(f_1,\ f_2)$ be the HVRFIF constructed in the Theorem 2 and $\omega_k+\tilde{\omega}_k<\frac{l_L}{\mathrm{L}_L},\ 1\le k\le n$. Then there exists positive constant $\mathrm{L}_1$, $\mathrm{L}_2$ and $\tau_1$, $\tau_2$ with $0<\tau_1,\tau_2\le 1$ such that for any $x_1,\ x_2\in I$,**

$$|f_1(x_1)-f_1(x_2)|\le\mathrm{L}_1|x_1-x_2|^{\tau_1},\ |f_2(x_1)-f_2(x_2)|\le\mathrm{L}_2|x_1-x_2|^{\tau_2}. \tag{5}$$

**Proof.** Let us denote $I_{r_1 r_2 \ldots r_m} = L_{r_m,\gamma(r_m)} \circ L_{r_{m-1},\gamma(r_{m-1})} \circ \cdots \circ L_{r_1,\gamma(r_1)}(\widetilde{I}_{\gamma(r_1)})$ . For any $x_1, x_2 \in I$ ( $0 \le x_1 < x_2 \le 1$), there exists $m$ such that $I_{r_1 r_2 \ldots r_m}$ is the largest interval contained in $[x_1, x_2]$. Then either $x_1 \in I_{r_2 \ldots r_m}$ or $x_2 \in I_{r_2 \ldots r_m}$.

Assume that $x_1 \in I_{r_2 \ldots r_m}$ . Then either $x_1 \in I_{s r_2 \ldots r_m}$, $s \le r_1 - 1$ , $x_2 \in I_{t r_2 \ldots r_m}$, $t \ge r_1 + 1$ or $x_1 \in I_{s r_2 \ldots r_m}$, $s \le r_1 - 1$, $x_2 \in I_{t r_2 + 1 \ldots r_m}$, $1 \le t \le n$.

Hence

$$
\begin{aligned}
| f_1(x_1) - f_1(x_2) | =& | (s_{r_m}(x_1) f_1(L^{-1}_{r_m,\gamma(r_m)}(x_1)) + s'_{r_m}(x_1) f_2(L^{-1}_{r_m,\gamma(r_m)}(x_1)) + q_{r_m,\gamma(r_m)}(L^{-1}_{r_m,\gamma(r_m)}(x_1))) - \\
& - (s_{r_m}(x_2) f_1(L^{-1}_{r_m,\gamma(r_m)}(x_2)) + s'_{r_m}(x_2) f_2(L^{-1}_{r_m,\gamma(r_m)}(x_2)) + q_{r_m,\gamma(r_m)}(L^{-1}_{r_m,\gamma(r_m)}(x_2))) | \\
& \le | s_{r_m}(x_1) f_1(L^{-1}_{r_m,\gamma(r_m)}(x_1)) - s_{r_m}(x_2) f_1(L^{-1}_{r_m,\gamma(r_m)}(x_2)) | + | s'_{r_m}(x_1) f_2(L^{-1}_{r_m,\gamma(r_m)}(x_1)) - \\
& - s'_{r_m}(x_2) f_2(L^{-1}_{r_m,\gamma(r_m)}(x_2)) | + | q_{r_m,\gamma(r_m)}(L^{-1}_{r_m,\gamma(r_m)}(x_1)) - q_{r_m,\gamma(r_m)}(L^{-1}_{r_m,\gamma(r_m)}(x_2)) |,
\end{aligned}
\tag{6}
$$

and

$$
\begin{aligned}
& | s_{r_m}(x_1) f_1(L^{-1}_{r_m,\gamma(r_m)}(x_1)) - s_{r_m}(x_2) f_1(L^{-1}_{r_m,\gamma(r_m)}(x_2)) | \\
& \le | s_{r_m}(x_1) f_1(L^{-1}_{r_m,\gamma(r_m)}(x_1)) - s_{r_m}(x_1) f_1(L^{-1}_{r_m,\gamma(r_m)}(x_2)) | + \\
& \quad + | s_{r_m}(x_1) f_1(L^{-1}_{r_m,\gamma(r_m)}(x_2)) - s_{r_m}(x_2) f_1(L^{-1}_{r_m,\gamma(r_m)}(x_2)) | \\
& \le | s_{r_m}(x_1) | \cdot | f_1(L^{-1}_{r_m,\gamma(r_m)}(x_1)) - f_1(L^{-1}_{r_m,\gamma(r_m)}(x_2)) | + | f_1(L^{-1}_{r_m,\gamma(r_m)}(x_2)) | \cdot | s_{r_m}(x_1) - s_{r_m}(x_2) | \\
& \le \bar{s}_{r_m} \, | f_1(L^{-1}_{r_m,\gamma(r_m)}(x_1)) - f_1(L^{-1}_{r_m,\gamma(r_m)}(x_2)) | + \overline{f_1} \cdot \mathrm{L}_{s_{r_m}} \, | x_1 - x_2 | \\
& \le \omega_{r_m} \, | f_1(L^{-1}_{r_m,\gamma(r_m)}(x_1)) - f_1(L^{-1}_{r_m,\gamma(r_m)}(x_2)) | + \overline{f_1} \cdot \mathrm{L}_{s_{r_m}} \, | x_1 - x_2 |,
\end{aligned}
$$

$$
\begin{aligned}
& | s'_{r_m}(x_1) f_2(L^{-1}_{r_m,\gamma(r_m)}(x_1)) - s'_{r_m}(x_2) f_2(L^{-1}_{r_m,\gamma(r_m)}(x_2)) | \\
& \le \omega_{r_m} \, | f_2(L^{-1}_{r_m,\gamma(r_m)}(x_1)) - f_2(L^{-1}_{r_m,\gamma(r_m)}(x_2)) | + \overline{f_2} \cdot \mathrm{L}_{s'_{r_m}} \, | x_1 - x_2 |,
\end{aligned}
$$

$$
\begin{aligned}
& | q_{r_m,\gamma(r_m)}(L^{-1}_{r_m,\gamma(r_m)}(x_1)) - q_{r_m,\gamma(r_m)}(L^{-1}_{r_m,\gamma(r_m)}(x_2)) | \\
& \le \mathrm{L}_{q_{r_m,\gamma(r_m)}} \cdot | L^{-1}_{r_m,\gamma(r_m)}(x_1) - L^{-1}_{r_m,\gamma(r_m)}(x_2) | = \mathrm{L}_{q_{r_m,\gamma(r_m)}} \cdot \frac{| \widetilde{I}_{\gamma(r_m)} |}{|I_{r_m}|} \cdot |x_1 - x_2| .
\end{aligned}
$$

Therefore, by (6), it follows that

$$
\begin{aligned}
| f_1(x_1) - f_1(x_2) | =& \left( \overline{f_1} \cdot \mathrm{L}_{s_{r_m}} + \overline{f_2} \cdot \mathrm{L}_{s'_{r_m}} + \mathrm{L}_{q_{r_m,\gamma(r_m)}} \cdot \frac{| \widetilde{I}_{\gamma(r_m)} |}{|I_{r_m}|} \right) | x_1 - x_2 | + \\
& + \omega_{r_m} ( | f_1(L^{-1}_{r_m,\gamma(r_m)}(x_1)) - f_1(L^{-1}_{r_m,\gamma(r_m)}(x_2)) | + | f_2(L^{-1}_{r_m,\gamma(r_m)}(x_1)) - f_2(L^{-1}_{r_m,\gamma(r_m)}(x_2)) | ) \\
& \le M_{r_m} \, | x_1 - x_2 | + \omega_{r_m} ( | f_1(L^{-1}_{r_m,\gamma(r_m)}(x_1)) - f_1(L^{-1}_{r_m,\gamma(r_m)}(x_2)) | + | f_2(L^{-1}_{r_m,\gamma(r_m)}(x_1)) - f_2(L^{-1}_{r_m,\gamma(r_m)}(x_2)) | ) ,
\end{aligned}
$$

where $M_k = \overline{f_1} \cdot \mathrm{L}_{s_k} + \overline{f_2} \cdot \mathrm{L}_{s'_k} + \mathrm{L}_{q_{k,\gamma(k)}} \cdot \dfrac{| \widetilde{I}_{\gamma(k)} |}{|I_k|}$.

Similarly, we have

$$
\begin{aligned}
| f_2(x_1) - f_2(x_2) | = \widetilde{M}_{r_m} \, | x_1 - x_2 | + \widetilde{\omega}_{r_m} ( & | f_1(L^{-1}_{r_m,\gamma(r_m)}(x_1)) - f_1(L^{-1}_{r_m,\gamma(r_m)}(x_2)) | + \\
& + | f_2(L^{-1}_{r_m,\gamma(r_m)}(x_1)) - f_2(L^{-1}_{r_m,\gamma(r_m)}(x_2)) | ) ,
\end{aligned}
$$

where $\widetilde{M}_k = \overline{f_1} \cdot \mathrm{L}_{\widetilde{s}_k} + \overline{f_2} \cdot \mathrm{L}_{\widetilde{s}'_k} + \mathrm{L}_{\widetilde{q}_{k,\gamma(k)}} \cdot \dfrac{| \widetilde{I}_{\gamma(k)} |}{|I_k|}$.

Hence by induction, we have

$$
\begin{aligned}
| f_1(x_1) - f_1(x_2) | &= M_{r_m} | x_1 - x_2 | + \omega_{r_m} \{ [ M_{r_{m-1}} | L^{-1}_{r_m,\gamma(r_m)}(x_1) - L^{-1}_{r_m,\gamma(r_m)}(x_2) | + \\
&\quad + \omega_{r_{m-1}} ( | f_1(L^{-1}_{r_{m-1},\gamma(r_{m-1})}(L^{-1}_{r_m,\gamma(r_m)}(x_1))) - f_1(L^{-1}_{r_{m-1},\gamma(r_{m-1})}(L^{-1}_{r_m,\gamma(r_m)}(x_2))) | + \\
&\quad + | f_2(L^{-1}_{r_{m-1},\gamma(r_{m-1})}(L^{-1}_{r_m,\gamma(r_m)}(x_1))) - f_2(L^{-1}_{r_{m-1},\gamma(r_{m-1})}(L^{-1}_{r_m,\gamma(r_m)}(x_2))) | ) ] + \\
&\quad + [ \tilde{M}_{r_{m-1}} | L^{-1}_{r_m,\gamma(r_m)}(x_1) - L^{-1}_{r_m,\gamma(r_m)}(x_2) | + \\
&\quad + \tilde{\omega}_{r_{m-1}} ( | f_1(L^{-1}_{r_{m-1},\gamma(r_{m-1})}(L^{-1}_{r_m,\gamma(r_m)}(x_1))) - f_1(L^{-1}_{r_{m-1},\gamma(r_{m-1})}(L^{-1}_{r_m,\gamma(r_m)}(x_2))) | + \\
&\quad + | f_2(L^{-1}_{r_{m-1},\gamma(r_{m-1})}(L^{-1}_{r_m,\gamma(r_m)}(x_1))) - f_2(L^{-1}_{r_{m-1},\gamma(r_{m-1})}(L^{-1}_{r_m,\gamma(r_m)}(x_2))) | ) ] \} \\
&= \left[ M_{r_m} + \frac{| \tilde{I}_{\gamma(r_m)} |}{| I_{r_m} |} \omega_{r_m} ( M_{r_{m-1}} + \tilde{M}_{r_{m-1}} ) \right] \cdot | x_1 - x_2 | + \\
&\quad + \omega_{r_m} ( \omega_{r_{m-1}} + \tilde{\omega}_{r_{m-1}} ) ( | f_1 \circ L^{-1}_{r_{m-1},\gamma(r_{m-1})} \circ L^{-1}_{r_m,\gamma(r_m)}(x_1) - f_1 \circ L^{-1}_{r_{m-1},\gamma(r_{m-1})} \circ L^{-1}_{r_m,\gamma(r_m)}(x_2) | + \\
&\quad + | f_2 \circ L^{-1}_{r_{m-1},\gamma(r_{m-1})} \circ L^{-1}_{r_m,\gamma(r_m)}(x_1) - f_2 \circ L^{-1}_{r_{m-1},\gamma(r_{m-1})} \circ L^{-1}_{r_m,\gamma(r_m)}(x_2) | ) \\
&\le \cdots \le \\
&\le \left[ M_{r_m} + \frac{\omega_{r_m}}{\omega_{r_m} + \tilde{\omega}_{r_m}} \sum_{j=4}^{m} \left( \prod_{k=j}^{m} \frac{| \tilde{I}_{\gamma(r_k)} |}{| I_{r_k} |} ( \omega_{r_k} + \tilde{\omega}_{r_k} ) \right) \cdot ( M_{r_{j-1}} + \tilde{M}_{r_{j-1}} ) \right] \cdot | x_1 - x_2 | + \\
&\quad + \frac{\omega_{r_m}}{\omega_{r_m} + \tilde{\omega}_{r_m}} \left( \prod_{j=3}^{m} ( \omega_{r_j} + \tilde{\omega}_{r_j} ) \right) ( | f_1 \circ L^{-1}_{r_3,\gamma(r_3)} \circ \cdots \circ L^{-1}_{r_m,\gamma(r_m)}(x_1) - f_1 \circ L^{-1}_{r_3,\gamma(r_3)} \circ \cdots \circ L^{-1}_{r_m,\gamma(r_m)}(x_2) | + \\
&\qquad + | f_2 \circ L^{-1}_{r_3,\gamma(r_3)} \circ \cdots \circ L^{-1}_{r_m,\gamma(r_m)}(x_1) - f_2 \circ L^{-1}_{r_3,\gamma(r_3)} \circ \cdots \circ L^{-1}_{r_m,\gamma(r_m)}(x_2) | ) \\
&\le \left[ M_{r_m} + \left( \sum_{j=4}^{m} ( M_{r_{j-1}} + \tilde{M}_{r_{j-1}} ) \, \delta^{m-j+1} \right) \right] \cdot | x_1 - x_2 | + \prod_{j=3}^{m} ( \omega_{r_j} + \tilde{\omega}_{r_j} ) ( 2\bar{f}_1 + 2\bar{f}_2 ) \\
&\le \left( \mathrm{M} \sum_{j=0}^{m-3} \delta^j \right) \cdot | x_1 - x_2 | + 2 \left( \prod_{j=3}^{m} \frac{| I_{r_j} |}{| \tilde{I}_{\gamma(r_j)} |} \right) \cdot ( \bar{f}_1 + \bar{f}_2 ) \delta^{m-2},
\end{aligned}
$$

where $\mathrm{M} = \max_k \{ M_k + \tilde{M}_k \}$, $\delta = \max_k \left\{ \frac{| \tilde{I}_{\gamma(k)} |}{| I_k |} ( \omega_k + \tilde{\omega}_k ) \right\}$.

Since $\prod_{j=1}^{m} \frac{| I_{r_j} |}{| \tilde{I}_{\gamma(r_j)} |} \cdot | \tilde{I}_{\gamma(r_1)} | \le | x_1 - x_2 | \le \prod_{j=3}^{m} \frac{| I_{r_j} |}{| \tilde{I}_{\gamma(r_j)} |} \cdot | \tilde{I}_{\gamma(r_3)} |$, we have

$$
\begin{aligned}
| f_1(x_1) - f_1(x_2) | &\le \left( \mathrm{M} \sum_{j=0}^{m-3} \delta^j \right) \cdot | x_1 - x_2 | + \frac{2 | \tilde{I}_{\gamma(r_2)} |}{| I_{r_1} | \cdot | I_{r_2} |} \cdot ( \bar{f}_1 + \bar{f}_2 ) \cdot \delta^{m-2} \cdot | x_1 - x_2 | \\
&\le \left( \mathrm{M} \sum_{j=0}^{m-3} \delta^j \right) \cdot | x_1 - x_2 | + \frac{2 | \tilde{I}_{\max} | \cdot ( \bar{f}_1 + \bar{f}_2 )}{| I_{\min} |^2} \cdot \delta^{m-2} \cdot | x_1 - x_2 | \\
&\le \mathrm{D} \left( \sum_{i=0}^{m-2} \delta^i \right) \cdot | x_1 - x_2 |, \qquad (7)
\end{aligned}
$$

where $\mathrm{D} = \max \left\{ \mathrm{M}, \ \frac{2 | \tilde{I}_{\max} | \cdot ( \bar{f}_1 + \bar{f}_2 )}{| I_{\min} |^2} \right\}$.

(1) If $\delta < 1$, then $\sum_{i=0}^{m-2} \delta^i \le \frac{1}{1-\delta}$. Then by (7), we have

$$
| f_1(x_1) - f_1(x_2) | \le \frac{\mathrm{D}}{1-\delta} | x_1 - x_2 |.
$$

Therefore, let us denote $\mathrm{L}_1 := \frac{\mathrm{D}}{1-\delta}$, $\tau_1 := 1$, then

$$| f_1(x_1) - f_1(x_2) | \le L_1 | x_1 - x_2 |^{\tau_1}$$

(2) If $\delta = 1$, then $\sum_{i=0}^{m-2} \delta^i = m-1$ and by (7), we have

$$| f_1(x_1) - f_1(x_2) | \le D(m-1) | x_1 - x_2 | .$$

Since $| x_1 - x_2 | \le \prod_{j=3}^{m} \frac{|I_{r_j}|}{|\tilde{I}_{\gamma(r_j)}|} \cdot | \tilde{I}_{\gamma(r_3)} | \le L_L^{m-2}$, we obtain $m - 2 \le \frac{\ln | x_1 - x_2 |}{\ln L_L}$.

Hence we have

$$| f_1(x_1) - f_1(x_2) | \le D\left( \frac{\ln | x_1 - x_2 |}{\ln L_L} + 1 \right) | x_1 - x_2 |$$

$$= D\left( | x_1 - x_2 | - \frac{| x_1 - x_2 |^{\alpha} \ln | x_1 - x_2 |}{| \ln L_L |} | x_1 - x_2 |^{1-\alpha} \right). \quad (8)$$

As $-| x_1 - x_2 |^{\alpha} \ln | x_1 - x_2 | \le \frac{1}{\alpha e}$ by Lemma1, we have

$$| f_1(x_1) - f_1(x_2) | \le D\left( | x_1 - x_2 | + \frac{1}{\alpha e | \ln L_L |} | x_1 - x_2 |^{1-\alpha} \right)$$

$$\le D\left( 1 + \frac{1}{\alpha e | \ln L_L |} \right) | x_1 - x_2 |^{1-\alpha} .$$

Let us denote $L_1 := D\left( 1 + \frac{1}{\alpha e | \ln L_L |} \right)$, $\tau_1 = 1 - \alpha$. Then

$$| f_1(x_1) - f_1(x_2) | \le L_1 | x_1 - x_2 |^{\tau_1}$$

(3) If $\delta > 1$, then $\sum_{i=0}^{m-2} \delta^i \le \frac{\delta^{m-2}}{1 - \frac{1}{\delta}} = \frac{\delta^{m-1}}{\delta - 1}$ and by (7), we get

$$| f_1(x_1) - f_1(x_2) | \le D \frac{\delta^{m-1}}{\delta - 1} | x_1 - x_2 | . \quad (9)$$

Since $| x_1 - x_2 | \le L_L^{m-2}$, we have $\ln | x_1 - x_2 | \le (m-2) \ln L_L$. Let us choose $\tau$ such that $0 < \tau \le \frac{\ln \delta}{\ln L_L} + 1$. Then

$$\frac{\ln(\delta^{m-1} | x_1 - x_2 |)}{\ln | x_1 - x_2 |} = \frac{(m-1) \ln \delta}{\ln | x_1 - x_2 |} + 1 \ge \frac{(m-1) \ln \delta}{(m-2) \ln L_L} + 1 \ge \frac{\ln \delta}{\ln L_L} + 1 \ge \tau .$$

Hence by (9), we get

$$| f_1(x_1) - f_1(x_2) | \le D \frac{1}{\delta - 1} | x_1 - x_2 |^{\tau} \le L_1 | x_1 - x_2 |^{\tau_1} ,$$

where $L_1 = \frac{D}{\delta - 1}$, $\tau_1 = \tau$.

Similarly, we can prove that there exists positive constant $L_2$ and $\tau_2 (0 < \tau_2 \le 1)$ such that for any $x_1, x_2 \in I$,

$$| f_2(x_1) - f_2(x_2) | \le L_2 | x_1 - x_2 |^{\tau_2} . \quad \square$$

## 4. Stability of the HVRFIF

In this section we consider the stability of the HVRFIF in the case of the interpolation points having small perturbations.

First, we prove the stability of the HVFIF according to perturbation of each coordinate of points in the data set $P$, and then all coordinates.

Let $I=[0,1]$. For $x_0=x_0^*<x_1^*<\cdots<x_n^*=x_n$, let us denote as follows:

$$I_i^*=[x_{i-1}^*,\,x_i^*],\ \ \tilde{I}_k^*=[x_{s(k)}^*,\,x_{e(k)}^*],\ \ L_{i,k}^*:\tilde{I}_k^*\to I_i^*.$$

Let us define the mapping $R:I\to I$ as follows:

$$R(x)=x_{i-1}^*+\frac{x_i^*-x_{i-1}^*}{x_i-x_{i-1}}(x-x_{i-1}),\ \text{for}\ \ x\in I_i. \tag{10}$$

Then $R(I_i)=I_i^*$, $R(\tilde{I}_k)=\tilde{I}_k^*$ and the following lemma holds.

**Lemma 2[17]. Let $R(x)$ be the mapping in (10). Then**

$$|R(x)-x|\le\max_i|x_i-x_i^*|\ \ \textbf{for}\ \ x\in I\,\textbf{.} \tag{11}$$

Now, let us denote as follows:

$$L_{i,k}^*=R\circ L_{i,k}\circ R^{-1},\ \ s_i^*=s_i\circ R^{-1},\ \ s_i'^*=s_i'\circ R^{-1},\ \ \tilde{s}_i^*=\tilde{s}_i\circ R^{-1},\ \ \tilde{s}_i'^*=\tilde{s}_i'\circ R^{-1}. \tag{12}$$

Let us denote by $\vec{f}^*=(f_1^*,\,f_2^*)$ the HVRFIF constructed using (12) and $q_{i,k}^*$, $\tilde{q}_{i,k}^*$ in the example 1 for a perturbed data set

$$P_{x^*}=\{(x_i^*,\ y_i,z_i)\in\mathrm{R}^3;i=0,1,\cdots,n\}(x_0=x_0^*<x_1^*<\cdots<x_n^*=x_n)\,.$$

Let us denote as follows:

$$\omega=\max_j\{\omega_j\},\ \tilde{\omega}=\max_j\{\tilde{\omega}_j\},\ |y|_{\max}=\max_i|y_i|,\ |z|_{\max}=\max_i|z_i|,\ |\tilde{I}|_{\min}=\min_{k=1,\cdots,l}\{|\tilde{I}_k|,\ |\tilde{I}_k^*|\}\,.$$

**Theorem 4. Let $f_1$, $f_1^*$ be HVRFIF for the data sets $P$ and $P_{x^*}$, respectively, and $\omega+\tilde{\omega}<\dfrac{l_L}{\mathrm{L}_L}$. Then there exists constant $N,\ \tau(0<\tau\le1)$ such that**

$$\|f_1-f_1^*\|_\infty\le\frac{\mathrm{L}_1(1+\omega-\tilde{\omega})+2\omega\mathrm{L}_2+\omega N}{1-\omega-\tilde{\omega}}\max_i|x_i-x_i^*|^\tau\,\textbf{.} \tag{13}$$

**Proof.** For any $x\in I$, there is $x_2\in\tilde{I}_l^*$ such that $x=L_{j,l}^*(x_2)$ and we get

$$R^{-1}(x)=R^{-1}(L_{j,l}^*(x_2))=(R^{-1}\circ L_{j,l}^*)(x_2)=(L_{j,l}\circ R^{-1})(x_2)\,.$$

Then, we have

$$\begin{aligned}f_1^*(x)&=s_j^*(x)f_1^*(L_{j,l}^{*-1}(x))+s_j'^*(x)f_2^*(L_{j,l}^{*-1}(x))+q_{j,l}^*(L_{j,l}^{*-1}(x))\\&=s_j\circ R^{-1}(x)f_1^*(x_2)+s_j'\circ R^{-1}(x)f_2^*(x_2)+q_{j,l}^*(x_2),\end{aligned} \tag{14}$$

$$|f_1(x)-f_1^*(x)|\le|f_1(x)-f_1\circ R^{-1}(x)|+|f_1\circ R^{-1}(x)-f_1^*(x)|\,. \tag{15}$$

By the smoothness of $f_1$ and Lemma 2, we get

$$|f_1(x)-f_1\circ R^{-1}(x)|\le\mathrm{L}_1|x-R^{-1}(x)|^{\tau_1}\le\mathrm{L}_1\max_i|x_i-x_i^*|^{\tau_1}\,, \tag{16}$$

$$| f_2(x) - f_2 \circ R^{-1}(x) | \le L_2 \max_i | x_i - x_i^* |^{\tau_2}.$$

In the other hand, we obtain

$$| (f_1 \circ R^{-1})(x) - f_1^*(x) | = | (s_j \circ R^{-1}(x) \cdot f_1 \circ L_{j,l}^{-1} \circ R^{-1}(x) + s'_j \circ R^{-1}(x) \cdot f_2 \circ L_{j,l}^{-1} \circ R^{-1}(x) +$$

$$+ q_{j,l} \circ L_{j,l}^{-1} \circ R^{-1}(x)) - (s_j \circ R^{-1}(x) \cdot f_1^*(x_2) + s'_j \circ R^{-1}(x) \cdot f_2^*(x_2) + q_{j,l}^*(x_2)) |$$

$$= | (s_j \circ R^{-1}(x) \cdot f_1 \circ R^{-1}(x_2) + s'_j \circ R^{-1}(x) \cdot f_2 \circ R^{-1}(x_2) + q_{j,l} \circ R^{-1}(x_2))$$

$$- (s_j \circ R^{-1}(x) \cdot f_1^*(x_2) + s'_j \circ R^{-1}(x) \cdot f_2^*(x_2) + q_{j,l}^*(x_2)) |$$

$$\le | (s_j \circ R^{-1}(x) | \cdot | f_1 \circ R^{-1}(x_2) - f_1^*(x_2) | + | s'_j \circ R^{-1}(x) | \cdot | f_2 \circ R^{-1}(x_2) - f_2^*(x_2) | +$$

$$+ | q_{j,l} \circ R^{-1}(x_2) - q_{j,l}^*(x_2) |$$

$$\le \bar{s}_j (| f_1 \circ R^{-1}(x_2) - f_1(x_2) | + | f_1(x_2) - f_1^*(x_2) |) + \bar{s}'_j (| f_2 \circ R^{-1}(x_2) - f_2(x_2) | + | f_2(x_2) - f_2^*(x_2) |) +$$

$$+ | q_{j,l} \circ R^{-1}(x_2) - q_{j,l}^*(x_2) |$$

$$\le \omega_j (L_1 | x_2 - R^{-1}(x_2) |^{\tau_1} + \| f_1 - f_1^* \|_\infty + L_2 | x_2 - R^{-1}(x_2) |^{\tau_2} + \| f_2 - f_2^* \|_\infty) + | q_{j,l} \circ R^{-1}(x_2) - q_{j,l}^*(x_2) |$$

$$\le \omega_j (L_1 \max_i | x_i - x_i^* |^{\tau_1} + \| f_1 - f_1^* \|_\infty + L_2 \max_i | x_i - x_i^* |^{\tau_2} + \| f_2 - f_2^* \|_\infty) + | q_{j,l} \circ R^{-1}(x_2) - q_{j,l}^*(x_2) | \quad (17)$$

Since $s_j^* \circ L_{j,l}^* = s_j \circ R^{-1} \circ R \circ L_{j,l} \circ R^{-1} = s_j \circ L_{j,l} \circ R^{-1}$, $s_j'^* \circ L_{j,l}^* = s'_j \circ L_{j,l} \circ R^{-1}$ we get

$$| q_{j,l} \circ R^{-1}(x_2) - q_{j,l}^*(x_2) | = | (-s_j \circ L_{j,l} \circ R^{-1}(x_2) \cdot g_{j,l} \circ R^{-1}(x_2) - s'_j \circ L_{j,l} \circ R^{-1}(x_2) \cdot g'_{j,l} \circ R^{-1}(x_2) +$$

$$+ h_j \circ L_{j,l} \circ R^{-1}(x_2)) - (-s_j^* \circ L_{j,l}^*(x_2) \cdot g_{j,l}^*(x_2) - s_j'^* \circ L_{j,l}^*(x_2) \cdot g_{j,l}'^*(x_2) + h_j^* \circ L_{j,l}^*(x_2)) |$$

$$\le | (-s_j \circ L_{j,l} \circ R^{-1}(x_2) \cdot g_{j,l} \circ R^{-1}(x_2) - s'_j \circ L_{j,l} \circ R^{-1}(x_2) \cdot g'_{j,l} \circ R^{-1}(x_2) + h_j \circ L_{j,l} \circ R^{-1}(x_2)) -$$

$$- (-s_j \circ L_{j,l} \circ R^{-1}(x_2) \cdot g_{j,l}^*(x_2) - s'_j \circ L_{j,l} \circ R^{-1}(x_2) \cdot g_{j,l}'^*(x_2) + h_j^* \circ L_{j,l}^*(x_2)) |$$

$$\le | s_j \circ L_{j,l} \circ R^{-1}(x_2) | \cdot | g_{j,l}^*(x_2) - g_{j,l} \circ R^{-1}(x_2) | + | s'_j \circ L_{j,l} \circ R^{-1}(x_2) | \cdot | g_{j,l}'^*(x_2) - g'_{j,l} \circ R^{-1}(x_2) | +$$

$$+ | h_j \circ L_{j,l} \circ R^{-1}(x_2) - h_j^* \circ L_{j,l}^*(x_2) |$$

$$\le \bar{s}_j \, | g_{j,l}^*(x_2) - g_{j,l} \circ R^{-1}(x_2) | + \bar{s}'_j \, | g_{j,l}'^*(x_2) - g'_{j,l} \circ R^{-1}(x_2) | + | h_j \circ L_{j,l} \circ R^{-1}(x_2) - h_j^* \circ L_{j,l}^*(x_2) |$$

$$\le \omega_j (| g_{j,l}^*(x_2) - g_{j,l} \circ R^{-1}(x_2) | + | g_{j,l}'^*(x_2) - g'_{j,l} \circ R^{-1}(x_2) |) +$$

$$+ | h_j \circ L_{j,l} \circ R^{-1}(x_2) - h_j^* \circ R \circ L_{j,l} \circ R^{-1}(x_2) |.$$

Because

$$| g_{j,l}^*(x_2) - g_{j,l} \circ R^{-1}(x_2) | =$$

$$= \left| \left( \frac{x_2 - x_{s(l)}^*}{x_{e(l)}^* - x_{s(l)}^*} y_{e(l)} + \frac{x_2 - x_{e(l)}^*}{x_{s(l)}^* - x_{e(l)}^*} y_{s(l)} \right) - \left( \frac{R^{-1}(x_2) - x_{s(l)}}{x_{e(l)} - x_{s(l)}} y_{e(l)} + \frac{R^{-1}(x_2) - x_{e(l)}}{x_{s(l)} - x_{e(l)}} y_{s(l)} \right) \right|$$

$$\le | y_{\max} | \cdot \left( \left| \frac{x_2 - x_{s(l)}^*}{x_{e(l)}^* - x_{s(l)}^*} - \frac{R^{-1}(x_2) - x_{s(l)}}{x_{e(l)} - x_{s(l)}} \right| + \left| \frac{x_2 - x_{e(l)}^*}{x_{s(l)}^* - x_{e(l)}^*} - \frac{R^{-1}(x_2) - x_{e(l)}}{x_{s(l)} - x_{e(l)}} \right| \right)$$

$$\le 2 | y_{\max} | \cdot \frac{(x_2 + x_{s(l)}) | x_{e(l)}^* - x_{e(l)} | + (x_2 + x_{e(l)}) | x_{s(l)}^* - x_{s(l)} | + (x_{s(l)}^* + x_{e(l)}^*) | x_2 - R^{-1}(x_2) |}{| x_{e(l)}^* - x_{s(l)}^* | \cdot | x_{e(l)} - x_{s(l)} |}$$

$$\le \frac{12 | y |_{\max}}{| \tilde{I} |_{\min}^2} \max_i | x_i - x_i^* |,$$

$$| g_{j,l}'^*(x_2) - g'_{j,l} \circ R^{-1}(x_2) | \le \frac{12 | z |_{\max}}{| \tilde{I} |_{\min}^2} \max_i | x_i - x_i^* |,$$

and

$$|h_j(x')-h_j^*\circ R(x')|=\left|\left(\frac{x'-x_{j-1}}{x_j-x_{j-1}}y_j+\frac{x'-x_j}{x_{j-1}-x_j}y_{j-1}\right)-\left(\frac{R(x')-x_{j-1}^*}{x_j^*-x_{j-1}^*}y_j-\frac{R(x')-x_j^*}{x_{j-1}^*-x_j^*}y_{j-1}\right)\right|=0,$$

we have

$$|q_{j,l}\circ R^{-1}(x_2)-q_{j,l}^*(x_2)|\le\frac{12\omega_j(|y|_{\max}+|z|_{\max})}{|\tilde{I}|_{\min}{}^2}\max_i|x_i-x_i^*|.$$

Hence by (15)-(17), we get

$$|f_1(x)-f_1^*(x)|\le \mathrm{L}_1\max_i|x_i-x_i^*|^{\tau_1}+\omega_j(\mathrm{L}_1\max_i|x_i-x_i^*|^{\tau_1}+\|f_1-f_1^*\|_\infty+\mathrm{L}_2\max_i|x_i-x_i^*|^{\tau_2}+$$

$$+\|f_2-f_2^*\|_\infty)+\frac{12\omega_j(|y|_{\max}+|z|_{\max})}{|\tilde{I}|_{\min}{}^2}\max_i|x_i-x_i^*|$$

$$\le(\mathrm{L}_1+\omega\mathrm{L}_1+\omega\mathrm{L}_2+\omega N)\max_i|x_i-x_i^*|^{\tau}+\omega\|f_1-f_1^*\|_\infty+\omega\|f_2-f_2^*\|_\infty,$$

where $\tau=\min\{\tau_1,\tau_2\}$, $N=\frac{12\omega_j(|y|_{\max}+|z|_{\max})}{|\tilde{I}|_{\min}{}^2}$. Therefore we have

$$\|f_1-f_1^*\|_\infty\le(\mathrm{L}_1+\omega\mathrm{L}_1+\omega\mathrm{L}_2+\omega N)\max_i|x_i-x_i^*|^{\tau}+\omega\|f_1-f_1^*\|_\infty+\omega\|f_2-f_2^*\|_\infty. \quad (18)$$

Similarly we can easily prove that

$$\|f_2-f_2^*\|_\infty\le(\mathrm{L}_2+\tilde{\omega}\mathrm{L}_1+\tilde{\omega}\mathrm{L}_2+\tilde{\omega}N)\max_i|x_i-x_i^*|^{\tau}+\tilde{\omega}\|f_1-f_1^*\|_\infty+\tilde{\omega}\|f_2-f_2^*\|_\infty,$$

and

$$\|f_2-f_2^*\|_\infty\le\frac{1}{1-\tilde{\omega}}\left((\mathrm{L}_2+\tilde{\omega}\mathrm{L}_1+\tilde{\omega}\mathrm{L}_2+\tilde{\omega}N)\max_i|x_i-x_i^*|^{\tau}+\tilde{\omega}\|f_1-f_1^*\|_\infty\right). \quad (19)$$

As a result, by (18) and (19), we have

$$\|f_1-f_1^*\|_\infty\le\frac{\mathrm{L}_1(1+\omega-\tilde{\omega})+2\omega\mathrm{L}_2+\omega N}{1-\omega-\tilde{\omega}}\max_i|x_i-x_i^*|^{\tau}.$$ □

Let us denote by $\vec{f}^*=(f_1^*,f_2^*)$ the HVRFIF constructed using $q_{i,k}^*$, $\tilde{q}_{i,k}^*$ in the example 1 for a perturbed data set $P_{y^*}=\{(x_i,\ y_i^*,\ z_i)\in\mathrm{R}^3; i=0,1,\cdots,n\}$.

**Theorem 5. Let $f_1$, $f_1^*$ be HVRFIF for the data sets $P$ and $P_{y^*}$, respectively, and $\omega+\tilde{\omega}<\frac{l_L}{\mathrm{L}_L}$. Then**

$$\|f_1-f_1^*\|_\infty\le\frac{1+2\omega-\tilde{\omega}}{1-\omega-\tilde{\omega}}\max_i|y_i-y_i^*|. \quad (20)$$

**Proof.** For any $x\in I$, there exists $x'\in\tilde{I}_k$ such that $x=L_{i,k}(x')$ and we have

$$f_1(x)=s_i(x)f_1(L_{i,k}^{-1}(x))+s_i'(x)f_2(L_{i,k}^{-1}(x))+q_{i,k}(L_{i,k}^{-1}(x))$$
$$=s_i(x)f_1(x')+s_i'(x)f_2(x')+q_{i,k}(x'),$$
$$f_1^*(x)=s_i(x)f_1^*(x')+s_i'(x)f_2^*(x')+q_{i,k}^*(x').$$

Then we get

$$|f_1(x)-f_1^*(x)|\le|s_i(x)|\cdot|f_1(x')-f_1^*(x')|+|s_i'(x)|\cdot|f_2(x')-f_2^*(x')|+|q_{i,k}(x')-q_{i,k}^*(x')|$$

$$\le\omega_i(\|f_1-f_1^*\|_\infty+\|f_2-f_2^*\|_\infty)+|q_{i,k}(x')-q_{i,k}^*(x')|. \quad (21)$$

In the other hand,

$$|q_{i,k}(x')-q^*_{i,k}(x')|=|(-s_i\circ L_{i,k}(x')\cdot g_{i,k}(x')-s'_i\circ L_{i,k}(x')\cdot g'_{i,k}(x')+h_i\circ L_{i,k}(x'))-$$
$$-(-s_i\circ L_{i,k}(x')\cdot g^*_{i,k}(x')-s'_i\circ L_{i,k}(x')\cdot g'_{i,k}(x')+h^*_i\circ L_{i,k}(x'))|$$
$$\le |s_i(x)|\cdot|g_{i,k}(x')-g^*_{i,k}(x')|+|h_i(x)-h^*_i(x)|$$
$$\le \omega_i\left|\frac{x-x_{s(k)}}{x_{e(k)}-x_{s(k)}}(y_{e(k)}-y^*_{e(k)})+\frac{x-x_{e(k)}}{x_{s(k)}-x_{e(k)}}(y_{s(k)}-y^*_{s(k)})\right|$$
$$+\left|\frac{x-x_{i-1}}{x_i-x_{i-1}}(y_i-y^*_i)+\frac{x-x_i}{x_{i-1}-x_i}(y_{i-1}-y^*_{i-1})\right|$$
$$\le \omega_i\left|\frac{x-x_{s(k)}}{x_{e(k)}-x_{s(k)}}+\frac{x-x_{e(k)}}{x_{s(k)}-x_{e(k)}}\right|\max_i|y_i-y^*_i|+\left|\frac{x-x_{i-1}}{x_i-x_{i-1}}+\frac{x-x_i}{x_{i-1}-x_i}\right|\max_i|y_i-y^*_i|$$
$$=(\omega_i+1)\max_i|y_i-y^*_i|. \quad (22)$$

Therefore, by (21) and (22) we have

$$\|f_1-f^*_1\|_\infty\le\omega(\|f_1-f^*_1\|_\infty+\|f_2-f^*_2\|_\infty)+(\omega+1)\max_i|y_i-y^*_i|. \quad (23)$$

Similarly, we get

$$\|f_2-f^*_2\|_\infty\le\tilde{\omega}(\|f_1-f^*_1\|_\infty+\|f_2-f^*_2\|_\infty)+(\tilde{\omega}+1)\max_i|y_i-y^*_i|,$$

and

$$\|f_2-f^*_2\|_\infty\le\frac{1}{1-\tilde{\omega}}\left(\tilde{\omega}\|f_1-f^*_1\|_\infty+(\tilde{\omega}+1)\max_i|y_i-y^*_i|\right). \quad (24)$$

Hence, by (23) and (24),

$$\|f_1-f^*_1\|_\infty\le\frac{\omega}{1-\tilde{\omega}}\|f_1-f^*_1\|_\infty+\frac{1+2\omega-\tilde{\omega}}{1-\tilde{\omega}}\max_i|y_i-y^*_i|,$$
$$\|f_1-f^*_1\|_\infty\le\frac{1+2\omega-\tilde{\omega}}{1-\omega-\tilde{\omega}}\max_i|y_i-y^*_i|.$$
□

Let us denote by $\vec{f}^*=(f^*_1,f^*_2)$ the HVRFIF constructed using $q^*_{i,k}$, $\tilde{q}^*_{i,k}$ in the example 1 for a perturbed data set $P_{z^*}=\{(x_i,\ y_i,\ z^*_i)\in\mathrm{R}^3;i=0,1,\cdots,n\}$. Then, we can similarly prove the following theorem.

**Theorem 6. Let $f_1$, $f^*_1$ be HVRFIF for the data sets $P$ and $P_{z^*}$, respectively, and $\omega+\tilde{\omega}<\frac{l_L}{\mathrm{L}_L}$. Then**

$$\|f_1-f^*_1\|_\infty\le\frac{\omega}{1-\omega-\tilde{\omega}}\max_i|z_i-z^*_i|. \quad (25)$$

Let us denote by $\vec{f}^*=(f^*_1,\ f^*_2)$ the HVRFIF constructed using the construction in section 2 for a perturbed data set $P^*=\{(x^*_i,\ y^*_i,\ z^*_i)\in\mathrm{R}^3;i=0,1,\cdots,n\}$. Then from the theorem 4 and theorem 6, we have the following theorem.

**Theorem 7. Let $f_1$, $f^*_1$ be HVRFIF for the data sets $P$ and $P^*$, respectively, and**

$\omega + \tilde{\omega} < \frac{l_L}{L_L}$. **Then**

$$\| f_1 - f_1^* \|_\infty \le \frac{[L_1(1+\omega-\tilde{\omega}) + 2\omega L_2 + \omega N]\max_i |x_i - x_i^*|^\tau + (1+2\omega-\tilde{\omega})\max_i |y_i - y_i^*| + \omega \max_i |z_i - z_i^*|}{1-\omega-\tilde{\omega}},$$

**where** $N$ **and** $\tau$ **are constants in theorem 4.**

## References

[1] C. H. Yun, Hidden variable recurrent fractal interpolation function with four function contractivity factors, arXiv: 1904.09110 [math.DS]

[2] M. F. Barnsley, Fractal function and interpolation, Constr. Approx. 2(1986) 303- 329.

[3] R. Malysz, The Minkowski dimension of the bivariate fractal interpolation surfaces, Chaos Solitons Fractals 27(2006) 1147–1156.

[4] W. Metzler, C.H. Yun, Construction of fractal interpolation surfaces on rectangular grids, Int. J. Bifur. Chaos. 20 (2010) 4079–4086.

[5] N. Zhao, Construction and application of fractal interpolation surfaces. Visual Comput 1996;12:132–46.

[6] G. Chen, The smoothness and dimension of fractal interpolation function, Appl. Math. J. Chinese Univ. Ser. B 11(4) (1996) 409–418.

[7] Z. Feng, H. Xie, On Stability of Fractal Interpolation, Fractals 6(3) (1998) 269-273.

[8] H.Y. Wang, On smoothness for a class of fractal interpolation surfaces, Fractals 14 (3) (2006) 223–230.

[9] H.Y. Wang, X.J. Li, Perturbation error analysis for fractal interpolation functions and their moments, Appl. Math.Lett. 21 (2008) 441–446.

[10] Z. Feng, Y. Feng, Z. Yuan, Fractal interpolation surfaces with function contractivity factors, Applied Mathematics Letters, 25(2012) 1896–1900.

[11] C. H. Yun and M. K. Ri, Hidden variable bivariate fractal interpolation functions and errors on perturbations of function contractivity factors, Asian-Eur.J.Math.(2017), doi:10.1142/s1793557119500219

[12] C.H. Yun, H.C. O, H.C. Choi, Construction of fractal surfaces by recurrent fractal interpolation curves, Chaos, Solitons & Fractals 66 (2014) 136–143.

[13] C. H. Yun, H. C. Choi and H. C. O, Construction of recurrent fractal interpolation surfaces with function contractivity factors and estimation of box-counting dimension on rectangular grids, *Fractals* **23**(3) (2015) 1550030.

[14] J. Ji, J. Peng, Analytical properties of bivariate fractal interpolation functions with contractivity factor functions, Internationa Journal of Computer Mathematics 90(3)(2013) 539-553.

[15] H.Y. Wang, Z.L. Fan, Analytical characteristics of fractal interpolation functions with function contractivity factors, Acta Math. Sinica (Chin. Ser.) 54 (1) (2011) 147–158 (in Chinese).

[16] H.Y. Wang, J.S. Yu, Fractal interpolation functions with variable parameters and their analytical properties, Journal of Approximation Theory 175 (2013) 1–18.

[17] C. H. Yun and M. K. Ri, Analytic properties of hidden variable bivariate fractal interpolation functions with four function contractivity factors, Fractals (2018) doi:10.1142/s0218348x1950018x

[18] M. F. Barnsley, J.H. Elton, D.P. Hardin, Recurrent iterated function systems. Constr Approx 1989;5:3–31 .

[19] A.E. Jacquin, Image coding based on a fractal theory of iterated contractive image transformations.' IEEE Trans. Image processing 1(1992) 18-30.

[20] P. Bouboulis, L. Dalla, V. Drakopoulos, Image compression using recurrent bivariate fractal interpolation surfaces, ' *Internat. J. Bifur. Chaos* 16(7)(2006): 2063-2071.

[21] W. Metzler, C. H. Yun and M. Barski, Image compression predicated on recurrent iterated function

systems, Essays Math. Stat. 2 (2012) 107−116

[22] P. Bouboulis, L. Dalla, V. Drakopoulos, Construction of recurrent bivariate fractal interpolation surfaces and computation of their box-counting dimension. J Approx Theory, 141(2006) 99–117.

[23] C.D. Lour, Fractal interpolation functions with partial self similarity, J. Math. Anal. Appl. 464 (2018) 911–923

[24] M. F. Barnsley, D. Hardin, P. Massopust, Hidden variable fractal interpolation functions, SIAM J. Math. Anal., 20(1989), 1218-1242.

[25] A.K.B. Chand, G.P. Kapoor. Smoothness analysis of coalescence hidden variable fractal interpolation functions. Int. J. of Non-Linear Sci. 3(2007) 15–26.

[26] A.K.B. Chand, G.P. Kapoor, Stability of affine coalescence hidden variable fractal interpolation functions, Nonlinear Anal. 68 (2008) 3757–3770.

[27] G.P. Kapoor, S. A. Prasad. Smoothness of hidden variable bivariate coalescence fractal interpolation surfaces. Int. J. of Bifurcat. Chaos. 19(7) (2009) 2321–2333.

[28] G.P.Kapoor, S. A. Prasad, Stability of Coalescence Hidden Variable Fractal Interpolation Surfaces, Int. J. of Non-Linear Sci. 9(3) (2010) 265-275.

[29] R. Uthayakumar and M. Rajkumar, Hidden variable bivariate fractal interpolation surfaces with function vertical scaling factor, International Journal of Pure and Applied Mathematics 106(5) (2016) 21-32.